\documentclass[11pt,a4paper]{article}
\usepackage{amscd,amssymb,amsthm, graphicx}
\begin{document}
\pagestyle{myheadings}
\def \e{\varepsilon}
\def \w{\widetilde}
\def \b{$\Box $}

\begin{center}
{\large \bf THE $~\e$ - REVISED SYSTEM OF THE RIGID BODY WITH THREE LINEAR CONTROLS}\\
\end{center}
\begin{center}
{\bf Dan COM\u ANESCU, Mihai IVAN and Gheorghe IVAN}
\end{center}

{\bf Abstract.} In this paper we introduce the $~\e$ - revised system associated to a Hamilton - Poisson system.
The $\e$ - revised system of the rigid body with three linear controls is defined and some of its geometrical and dynamical properties are investigated.
{\footnote{
Key words and phrases: almost metriplectic system, $~\e$- revised system, rigid body}

\section{ Introduction}

\hspace*{0.7cm}It is well known that many dynamical systems can be formulated using a Poisson structure (see for instance, R. Abraham and J. E.
Marsden [1] and  M. Puta [11]).

The metriplectic systems was introduced by P. J. Morrison in the paper [8]. These systems combine both the conservative and dissipative systems.

A metriplectic system is a differential system of the form
$~\dot{x} = P dH + g dC,~ $ where $ P $ is a Poisson tensor on a manifold $ M,~{\bf g} $ is a symmetric tensor of type $ ( 2,0 ) $ on $ M, $ and $ H $ and $ C $ are two smooth functions on $ M $ with the additional requirements:\\
$(a)~~~ P dC = 0 ; ~~~ (b) ~~~ g dH = 0 ~~~\hbox{and}~~~ (c) ~~~ dC\cdot g dC \leq 0.$

The differential systems of the form $~\dot{x} = P dH + g dC~$ which satisfies only the conditions (a) and (b) are called almost metriplectic
systems ( see  Fish, [2]; Marsden, [7]; Ortega and Planas - Bielsa, [9] ). An interesting class of almost metriplectic systems are so-called the
revised dynamical systems associated to Hamilton-Poisson systems (see  Gh. Ivan and D. Opri\c s, [5]).

The control of the rotation rigid body is one of the problems with a large practical applicability. For this reason, in this paper we study
the $ \e $ - revised dynamical system associated to the rigid body with three linear controls.

\section { Almost metriplectic systems}

\hspace*{0.7cm}We start this section with the presentation of the concept of almost metriplectic manifold (see Ortega and Planas- Bielsa, [9]).

Let $ M $ be a smooth manifold of dimension $ n $ and let $ C^{\infty}(M) $ be the ring of smooth real-valued functions on $ M $.

A {\it Leibniz manifold} is a pair $ ( M , [\cdot, \cdot ] ), $ where  $ [\cdot,\cdot ] $ is a Leibniz bracket on $ M $, that is
$ [\cdot,\cdot ] : C^{\infty}(M)\times C^{\infty}(M) \to C^{\infty}(M) $
is a {\bf R} - bilinear operation satisfying the following two conditions:\\
$(i)~~~~$ the {\it left Leibniz rule}:\\[-0.4cm]
$$[ f_{1}\cdot f_{2} , f_{3} ] = [ f_{1} , f_{3} ]\cdot f_{2} + f_{1} \cdot [ f_{2} , f_{3} ] ~~~\hbox{for all}~~~ f_{1}, f_{2}, f_{3} \in C^{\infty}(M);$$\\[-0.4cm]
$(ii)~~~$ the {\it right Leibniz rule}:\\[-0.4cm]
$$[ f_{1}, f_{2} \cdot f_{3} ] = [ f_{1} , f_{2} ]\cdot f_{3} + f_{2} \cdot [ f_{1} , f_{3} ] ~~~\hbox{for all}~~~ f_{3}, f_{3}, f_{3} \in C^{\infty}(M)$$\\[-0.4cm]
where "$ \cdot $" denotes the ordinary multiplication of functions.

\noindent Let $~P~$ and $~{\bf g } $ be two tensor fields of type $ (2,0) $ on $ M $ and $ \e\in {\bf R} $ be a parameter.

 \noindent We define the map $ [\cdot,(\cdot,\cdot) ]_{\e} : C^{\infty}(M)\times (C^{\infty}(M)\times C^{\infty}(M)) \to C^{\infty}(M) $  by:
\begin{equation}
[ f,( h_{1},h_{2})]_{\e} = P( d f, d h_{1} ) + \e {\bf g } ( d f, d h_{2} ) , ~~\hbox{ for all }~~ f,h_{1},h_{2} \in C^{\infty}(M). \label{1}
\end{equation}
\smallskip

{\bf Proposition 2.1.} {\it The map $ [\cdot,(\cdot,\cdot) ]_{\e} $  given by $ (1) $ satisfy the following relations:\\
$(i)~~~~~[ a f_{1} + b f_{2}, ( h_{1}, h_{2} )]_{\e} = a [ f_{1} , ( h_{1}, h_{2} )]_{\e} + b [ f_{2}, ( h_{1}, h_{2} )]_{\e};$\\[0.2cm]
$(ii)~~~~[ f , a ( h_{1}, h_{2} ) + b ( h_{1}^{\prime}, h_{2}^{\prime})]_{\e} = a [ f, ( h_{1}, h_{2} )]_{\e} + b [ f, ( h_{1}^{\prime}, h_{2}^{\prime} )]_{\e};$\\[0.2cm]
$(iii)~~~[ ff_{1} , ( h_{1}, h_{2} )]_{\e} = f [ f_{1}, ( h_{1}, h_{2} )]_{\e} + f_{1} [ f, ( h_{1}, h_{2})]_{\e};$\\[0.2cm]
$(iv)~~~~[ f , h ( h_{1}, h_{2} )]_{\e} = h [ f , ( h_{1}, h_{2} )]_{\e} + h_{1} P( d f, d h ) + \e  h_{2} {\bf g } ( d f, d h),$\\[0.2cm]
for all  $ f, f_{1} , f_{2}, h_{1}, h_{2}, h_{1}^{\prime}, h_{2}^{\prime}\in C^{\infty}(M) $ and $ a,b\in {\bf R} $.}
\bigskip

{\bf Proof.} Applying the properties of the differential of functions and using that $ P $ and $ {\bf g} $ are $ {\bf R} $- bilinear maps, it is easy to establish the relations
$(i) - (iv). $\hfill\b
\bigskip

 We consider the map $ [[\cdot,\cdot ]]_{\e} : C^{\infty}(M)\times C^{\infty}(M)\to C^{\infty}(M) $ defined by:
\begin{equation}
[[ f, h ]]_{\e} = [ f, ( h, h )]_{\e}, ~~\hbox{for all}~~ f,h \in C^{\infty}(M).\label{2}
\end{equation}

Therefore, the map $ [[\cdot,\cdot ]]_{\e} : C^{\infty}(M)\times C^{\infty}(M)\to C^{\infty}(M) $ is given by:
\begin{equation}
[[ f, h ]]_{\e} =P( d f, d h ) + \e  {\bf g } ( d f, d h ), ~~\hbox{for all}~~ f,h \in C^{\infty}(M).\label{3}
\end{equation}
\smallskip

{\bf Proposition 2.2.} {\it The bracket $ [[\cdot,\cdot ]]_{\e} $ on $ M $ given by $(3)$ verify the right Leibniz rule:\\[-0.2cm]
$$[[ f , h h^{\prime}]]_{\e} = h [[ f , h^{\prime} ]]_{\e} + h^{\prime} [[  f,  h ]]_{\e},~~~\hbox{for all}~~~ f, h, h^{\prime}\in C^{\infty}(M). $$}\\[-0.8cm]
\bigskip

{\bf Proof.} Indeed, $ [[ f , h h^{\prime}]]_{\e} = [ f, ( hh^{\prime}, hh^{\prime})]_{\e} = [ f, h(h^{\prime}, h^{\prime})]_{\e}. $
Putting $ h_{1} = h_{2} = h^{\prime} $ in the relation $(iv)$ from Proposition 2.1, we have succesive:\\
$ [ f, h(h^{\prime}, h^{\prime})]_{\e} = h [ f, (h^{\prime}, h^{\prime})] + h^{\prime} P(df,dh) + \e  h^{\prime} {\bf g}(df,dh) =
h [ f, (h^{\prime}, h^{\prime})]_{\e} + h^{\prime} ( P(df,dh) + \e {\bf g}(df,dh)) = h [ f, (h^{\prime}, h^{\prime})]_{\e} + h^{\prime} [ f, (h,h)]_{\e} = h [[ f, h^{\prime}]]_{\e} + h^{\prime} [[ f, h ]]_{\e}.$\hfill\b
\bigskip

By Proposition 2.1. (i), (ii) and (iii) and Proposition 2.2, we have that the map $ [[\cdot, \cdot]]_{\e} $ given by $ (3) $ is a Leibniz bracket on $ M.$
Hence, $[[\cdot, \cdot ]]_{\e} $ defines a Leibniz structure on the manifold $ M $ and  $ (M, P, {\bf g}, [[\cdot, \cdot ]]_{\e} ) $ is a Leibniz manifold for each $ \e \in {\bf R}$.

A Leibniz manifold $ (M, P,{\bf g }, [[ \cdot, \cdot ]]_{\e}) $ such that $ P $ is a skewsymmetric tensor field and $ {\bf g} $   is a symmetric tensor field is called {\it almost metriplectic manifold}. In other words,
given a skewsymmetric tensor field $ P $ of type $(2,0)$ and a symmetric tensor field ${\bf g} $  of type $(2,0)$ on a manifold $ M $, we can define an almost metriplectic structure on $ M $.

If the tensor field $ P $ is Poisson and the tensor field $ {\bf g } $ is nondegenerate, then  $ (M, P ,{\bf g }, [[ \cdot, \cdot ]]_{\e}) $ is
a {\it metriplectic manifold}, see Ortega $\&$ Planas - Bielsa [9].
\bigskip

{\bf Proposition 2.3.} {\it Let $ (M, P,{\bf g }, [[ \cdot, \cdot ]]_{\e}) $ be an almost metriplectic manifold. If there exist
 $ h_{1}, h_{2} \in C^{\infty}(M) $ such that $ P( df, dh_{2} ) = 0 $ and $ {\bf g}( df, dh_{1} ) = 0 $ for all $ f\in C^{\infty}(M) $, then the bracket $ [[\cdot, \cdot ]] $  given by $ (3) $ satisfies the relation:
\begin{equation}
[[ f, h_{1}+ h_{2} ]]_{\e} = [ f, ( h_{1} , h_{2} ) ]_{\e} , ~~\hbox{ for all}~~ f\in C^{\infty}(M). \label{4}
\end{equation}}\\[-0.8cm]
\smallskip

{\bf Proof.} Indeed, $~ [[ f, h_{1}+ h_{2} ]]_{\e} = P( df, d(h_{1} + h_{2}) ) + \e {\bf g}( df, d(h_{1} + h_{2}) ) =
P( df, d h_{1} + d h_{2} ) + \e {\bf g}( df, d h_{1} + d h_{2} ) = P( df, d h_{1}) + P( df, d h_{2} ) + \e {\bf g}( df, d h_{1} ) + \e {\bf g}( df, d h_{2} ) =
P( df, d h_{1}) + \e  {\bf g}( df, d h_{2} ) = [ f, ( h_{1} , h_{2} ) ]_{\e}.$\hfill\b
\bigskip

Let $ (M, P,{\bf g }, [[ \cdot, \cdot ]]_{\e}) $ be an almost metriplectic manifold and let
 $ h_{1}, h_{2} \in C^{\infty}(M) $ two functions such that $ P( df, dh_{2} ) = 0 $ and $ {\bf g}( df, dh_{1} ) = 0 $ for all $ f\in C^{\infty}(M). $ The
 vector field $ X_{h_{1}h_{2}} $ given by:\\[-0.4cm]
$$X_{h_{1} h_{2}}(f) = [[ f,  h_{1}+ h_{2} ]]_{\e} ~~\hbox{ for any }~~f\in C^{\infty}(M)$$\\[-0.4cm]
is called the {\it Leibniz vector field } associated to the triple $ (h_{1}, h_{2}, \e )$ on $ M.$

Taking account into Proposition 2.3 and $(1),~ X_{h_{1} h_{2}} $ is given by:
\begin{equation}
X_{h_{1} h_{2}}(f)= [ f, ( h_{1} , h_{2} ) ]_{\e} =  P( d f, d h_{1} ) + \e {\bf g } ( d f, d h_{2} ) , ~~\hbox{ for all}~~ f\in C^{\infty}(M). \label{5}
\end{equation}

In local coordinates on $ M, $ the differential system given by:
\begin{equation}
{\dot x}^{i} = [[ x^{i}, h_{1} + h_{2}]]_{\e} = [ x^{i} , ( h_{1}, h_{2} ) ]_{\e} \label{6}
\end{equation}
where
\begin{equation}
[ x^{i} , ( h_{1} , h_{2} ) ]_{\e} = X_{h_{1}h_{2}}(x^{i}) = P^{ij} \frac {\partial h_{1}}{\partial x^{j}} + \e g^{ij}\frac{\partial h_{2}}{\partial x^{j}}, ~ i,j=\overline{1,n} \label{7}
\end{equation}
with $ P^{ij} = P(dx^{i}, dx^{j}) $ and $ g^{ij} = {\bf g}(dx^{i}, dx^{j}),$  is called the {\it almost metriplectic system } on $ M $ associated to the Leibniz vector field $ X_{h_{1}h_{2}} $
with the bracket $ [[\cdot, \cdot ]]_{\e}.$

 We denote the matrix of the tensor fields $ P $ and ${\bf g}$  respectively by $ P=( P^{ij} ) $ and $ g = (g^{ij}).$ We have that $ P $ is a skewsymmetric matrix and $ g $ is a symmetric matrix.

We give now a way for to produce almost metriplectic manifolds.
\bigskip

{\bf Proposition 2.4.} {\it For a skewsymmetric tensor $ P $ of type $ ( 2,0 ) $ on a manifold $ M $ and two functions  $ h_{1}, h_{2} \in C^{\infty}(M) $ such that $ P( df, dh_{2} ) = 0 $ for all $ f\in C^{\infty}(M) $,
there exists a symmetric tensor $ {\bf g } $ of type $ ( 2,0 ) $ on $ M $ such that $ {\bf g}( df, dh_{1} ) = 0 $ for all $ f\in C^{\infty}(M) $ and $ ( M , P, {\bf g},[[\cdot, \cdot ]]_{\e} ) $  is an almost metriplectic manifold.}

{\bf Proof.} In a system of local coordinates on $ M ,$ let $ g = ( g^{ij} ) $ the matrix of the symmetric tensor $ {\bf g} $ which must to be determined.
 Then, the components $ g^{ij},~ i,j=\overline{1.n} $ verify the system of differential equations
$~ g^{ij}\displaystyle\frac{\partial h_{1}}{\partial x^{j}} = 0 , ~i,j=\overline{1,n}.~$

 In a chart $ U $ such that $~ \displaystyle\frac{\partial h_{1}}{\partial x^{j}}(x)\neq 0~,$  the components $ g^{ij} $ are given by:
\begin{equation}
\left\{\begin{array}{ccl}
g^{ii}(x) & = & - \sum\limits_{k=1,~k\neq i}^{n} (\displaystyle\frac{\partial h_{1}}{\partial x^{k}})^{2}\cr
g^{ij}(x) & = & \displaystyle\frac{\partial h_{1}}{\partial x^{i}}\displaystyle\frac{\partial h_{1}}{\partial x^{j}},~~~\hbox{for}~~i\neq j\cr
\end{array}\right.\label{8}
\end{equation}

Applying now Proposition 2.3 we obtain the result.\hfill\b
\bigskip

Proposition 2.4 is useful when we consider the $\e$- revised system of a Hamilton-Poisson system.

For this, let be a Hamilton-Poisson system on $ M $ described by the Poisson tensor  $ P $ having the matrix $ P = (P^{ij} ) $ and by
the Hamiltonian function $ h_{1}\in C^{\infty}(M) $  with the Casimir function $ h_{2}\in C^{\infty}(M) $ ( i.e. $ P^{ij} \frac {\partial h_{2}}{\partial x^{j}}= 0 $ for $~i,j=\overline{1,n} $ ).
The differential equations of the Hamilton-Poisson system are the following:
\begin{equation}
{\dot x}^{i} = P^{ij} \frac {\partial h_{1}}{\partial x^{j}}, ~~i,j=\overline{1,n}. \label{9}
\end{equation}

Using $ (8) ,$ we determine the matrix $ g = ( g^{ij} ) $ and we have:
\begin{equation}
g^{ij} \frac {\partial h_{1}}{\partial x^{j}}= 0, ~~i,j=\overline{1,n}. \label{10}
\end{equation}

Applying now  Proposition 2.4, for each $ \e\in {\bf R}, $ we obtain an almost metriplectic structure on $ M $ associated to system $ (9) .$ The differential system
associated to this structure is called the {\it $ \e $ - revised system} of the Hamilton - Poisson system.

Hence, the $ \e $ - revised system of the Hamilton - Poisson system defined by $ (9) $ is:
\begin{equation}
{\dot x}^{i} = P^{ij} \frac {\partial h_{1}}{\partial x^{j}} +  \e g^{ij} \frac {\partial h_{2}}{\partial x^{j}} , ~~i,j=\overline{1,n}. \label{11}
\end{equation}
The terms $ g^{ij} \frac {\partial h_{2}}{\partial x^{j}},~~i,j=\overline{1,n}~ $ from the $ \e $ - revised system $ (11) $ describe a {\it cube perturbation} of the Hamilton - Poisson system.
\bigskip

{\bf Remark 2.1.} We observe that the $ 0 $- revised system $ (11) $ coincide with the Hamilton - Poisson system $(9)$.\hfill\b

\section {The $ \e $ - revised system associated to the rigid body with three linear controls}

\hspace*{0.7cm}The rigid body equations with three linear controls ( see, M. Puta and D. Com\u anescu [12] ) are given by:
\begin{equation}
\left\{ \begin{array}{ccl}
{\dot x}^{1} & = & (a_{3}- a_{2}) x^{2}x^{3} + cx^{2} - b x^{3}\cr
{\dot x}^{2} & = & ( a_{1} - a_{3} ) x^{1} x^{3} - c x^{1} + a x^{3} \cr
{\dot x}^{3} & = & ( a_{2} - a_{1}) x^{1}x^{2} + b x^{1} - a x^{2} \cr
\end{array}\right.\label{12}
\end{equation}
where $ x(t)= ( x^{1}(t), x^{2}(t), x^{3}(t) )\in {\bf R}^{3} $ and
$ a_{1} =\displaystyle\frac{1}{I_{1}} , ~
a_{2} =\displaystyle\frac{1}{I_{2}}, ~ a_{3} =\displaystyle\frac{1}{I_{3}} $ with
$ I_{1} > I_{2} > I_{3} > 0 $ ( $~ I_{1}, I_{2}, I_{3} $ being the principal moments of inertia of the body )
and $ a, b, c \in {\bf R} $ are feedback parameters. We have $ 0 < a_{1} < a_{2} < a_{3}. $

The dynamics $ (12) $ is described by the Poisson tensor $ \Pi $ and by the Hamiltonian $ H $ on $ {\bf R}^{3} $ given by:
\begin{equation}
\Pi(x) = \left ( \begin{array}{ccc}
0 & - x^{3} & x^{2}\\
x^{3} & 0 & - x^{1}\\
- x^{2} & x^{1} & 0\\
\end{array}\right ),\label{13}
\end{equation}
\begin{equation}
H(x) = \displaystyle\frac{1}{2} [ a_{1}(x^{1})^{2} + a_{2}(x^{2})^{2} + a_{3}(x^{3})^{2} ] + a x^{1} + b x^{2} + c x^{3}.\label{14}
\end{equation}

Using $(13)$ and $(14)$, the dynamics $(12)$ can be written in the matrix form:
\begin{equation}
{\dot x}(t) = \Pi(x(t))\cdot \nabla H(x(t)), \label{15}
\end{equation}
where $ {\dot x}(t)= ( {\dot x}^{1}(t), {\dot x}^{2}(t),{\dot x}^{3}(t) )^{T} $ and $ \nabla H(x(t)) $ is the gradient of the Hamiltonian function $ H $ with respect
to the canonical metric on $ {\bf R}^{3} $.

Therefore, the dynamics $(12)$ has the Hamilton-Poisson formulation $~( {\bf R}^{3}, \Pi, H ),~$ where $ \Pi $ and $ H $ are given
by $ (13) $ and $ (14). $

The function $ C\in C^{\infty}({\bf R}^{3}) $ given by:
\begin{equation}
C(x) = \displaystyle\frac{1}{2}[ (x^{1})^{2} + ( x^{2})^{2} + (x^{3})^{2}] \label{16}
\end{equation}
is a Casimir of the configuration $~( {\bf R}^{3}, \Pi ),~$ i.e.
\begin{equation}
C(x)\cdot \nabla H(x)= O. \label{17}
\end{equation}

Applying the relations $ (8) $ for $ P = \Pi, h_{1}(x) = H(x) $ and $ h_{2}(x)= C(x), $ the symmetric tensor $ {\bf g} $ is given by the matrix:\\
$$g =\left(\begin{array}{ccc}
- (a_{2}x^{2} + b)^{2} - (a_{3}x^{3} + c)^{2} & (a_{1}x^{1} + a)(a_{2} x^{2} + b) & (a_{1}x^{1} + a)(a_{3} x^{3}+ c) \\
(a_{1}x^{1} + a)(a_{2} x^{2} + b)   & - (a_{1}x^{1} + a)^{2} - (a_{3}x^{3} + c)^{2}  & (a_{2}x^{2} + b)(a_{3} x^{3} + c) \\
(a_{1}x^{1} + a)(a_{3} x^{3}+ c)   &  (a_{2}x^{2} + b)(a_{3} x^{3} + c)  & - (a_{1}x^{1} + a)^{2} - (a_{2}x^{2} + b)^{2}\\
\end{array} \right)$$\\

since $~ \displaystyle\frac{\partial h_{1}}{\partial x^{1}} = a_{1}x^{1} + a,~~~ \displaystyle\frac{\partial h_{1}}{\partial x^{2}} = a_{2} x^{2} + b,~~~
\displaystyle\frac{\partial h_{1}}{\partial x^{3}} = a_{3}x^{3} + c.$\\

We have
\begin{equation}
g(x)\cdot \nabla h_{2}(x) = ( v_{1}(x), v_{2}(x), v_{3}(x) )^{T} \label{18}
\end{equation}
where
\begin{equation}
\left\{\begin{array}{l}
v_{1}(x)  = -[(a_{2}x^{2} + b)^{2} + (a_{3}x^{3} + c)^{2}]x^{1} + (a_{1}x^{1} + a)[(a_{2}x^{2} + b)x^{2} + (a_{3}x^{3} + c)x^{3}] \\

v_{2}(x)  = -[(a_{1}x^{1} + a)^{2} + (a_{3}x^{3} + c)^{2}]x^{2} + (a_{2}x^{2} + b)[(a_{1}x^{1} + a)x^{1} + (a_{3}x^{3} + c)x^{3}] \\

v_{3}(x)  = -[(a_{1}x^{1} + a)^{2} + (a_{2}x^{2} + b)^{2}]x^{3} + (a_{3}x^{3} + c)[(a_{1}x^{1} + a)x^{1} + (a_{2}x^{2} + b)x^{2}] \\
\end{array}\right.\label{19}
\end{equation}

The $ \e $ - revised system associated to dynamics $(12)$ is:
\begin{equation}
\left\{ \begin{array}{l}
{\dot x}^{1} = [(a_{3} - a_{2}) x^{2} x^{3} + c x^{2} - b x^{3}] + \e v_{1}(x) \\

{\dot x}^{2} = [(a_{1} - a_{3}) x^{1} x^{3} - c x^{1} + a x^{3}] + \e v_{2}(x) \\

{\dot x}^{3} = [(a_{2} - a_{1}) x^{1} x^{2} + b x^{1} - a x^{2}] + \e v_{3}(x)\\

\end{array}\right.\label{20}
\end{equation}

The differential system $ (20) $ is called the {\it $ \e $ - revised system } of
the rigid body with three linear controls.
Taking  $ a = b = c = 0 $ and $ \e = 1 $ in $(20)$,  we obtain the revised system of the free rigid body, see [5].\\

{\bf Vector writing of the $\e$- revised system $(20)$}.
We introduce the following notations:\\[0.1cm]
${\bf x} = ( x^{1}, x^{2}, x^{3} ),~~{\bf v} = ( v_{1}, v_{2}, v_{3} ),~~{\bf a} = ( a, b, c ),~~
{\bf m} ({\bf x}) = ( a_{1} x^{1} + a, a_{2} x^{2} + b, a_{3} x^{3} + c ).$\\[0.1cm]

For all $ {\bf u} = ( u_{1}, u_{2}, u_{3} ), {\bf w} = ( w_{1}, w_{2}, w_{3} )\in {\bf R}^{3}, $ the following relation holds:
\begin{equation}
{\bf u}\cdot [ {\bf w}\times ( {\bf u}\times {\bf w})] = ( {\bf u}\times {\bf w})^{2},~~~\hbox{with}~~~( {\bf u}\times {\bf w})^{2}=( {\bf u}\times {\bf w})\cdot ( {\bf u}\times {\bf w})\label{21}
\end{equation}
where $"\times"$ and $"\cdot"$ denote the cross product resp. inner product in $ {\bf R}^{3}; $ that is:\\[0.1cm]
${\bf u}\times {\bf w} = ( u_{2}w_{3} - u_{3}w_{2}, u_{3}w_{1} - u_{1}w_{3},u_{1}w_{2} - u_{2}w_{1}), ~~~ {\bf u}\cdot {\bf w} = u_{1}w_{1}+ u_{2}w_{2} + u_{3}w_{3}.$\\[0.1cm]

With the above notations, the dynamics $ (12) $ has the vector form:
\begin{equation}
\dot{{\bf x}} = {\bf x} \times {\bf m} ({\bf x}) .\label{22}
\end{equation}

It is not hard to verify the following equality:
\begin{equation}
{\bf v} =( {\bf x}\times {\bf m}({\bf x}))\times {\bf m}({\bf x}). \label{23}
\end{equation}

Using the relations $ (22), (23) $ and $(20)$, we can written the $\e$- revised system $ (20) $ in the vector form:
\begin{equation}
\dot{{\bf x}} = {\bf x} \times {\bf m} ({\bf x}) +\e [({\bf x}\times {\bf m}({\bf x}))\times {\bf m}({\bf x})]. \label{24}
\end{equation}

\section {The equilibrium points of the $ \e $ - revised system}

\hspace*{0.7cm}The equilibrium points of the Hamilton - Poisson system $ (12) $ ( or $(22)$ ) are solutions of the vector equation:
\begin{equation}
{\bf x} \times {\bf m} ({\bf x}) = {\bf 0} .\label{25}
\end{equation}

The equilibrium points of the $ \e $ - revised system $(20)$ ( or $(24)$ ) are solutions of the vector equation:
\begin{equation}
{\bf x} \times {\bf m} ({\bf x}) +\e [({\bf x}\times {\bf m}({\bf x}))\times {\bf m}({\bf x})] = {\bf 0}. \label{26}
\end{equation}
\smallskip

{\bf Theorem 4.1.} {\it The Hamilton - Poisson system $ (12) $ and its revised system $ (20) $ have the same equilibrium points}.
\bigskip

{\bf Proof.} Let $ {\bf x}_{0} $ be an equilibrium point of the system $ (12).$ According with $(25)$ follows ${\bf x}_{0} \times {\bf m}({\bf
x}_{0}) = {\bf 0}. $ We have that $ {\bf x}_{0} $ is a solution of the vector equation $(26)$, since ${\bf }x_{0} \times {\bf m}({\bf x}_{0}) +
({\bf x}_{0}\times {\bf m}({\bf x}_{0}))\times {\bf m}({\bf x}_{0}) = {\bf 0} + \textbf{0}\times {\bf m}({\bf x}_{0}) = {\bf 0}. $ Hence $ {\bf
x}_{0} $ is an equilibrium point of the $\e $- revised system $ (20) $.

Conversely, let $ {\bf x}_{0} $ be an equilibrium point for $ (20) $. Using $(25)$ it follows

\noindent $(a)~~~~~{\bf x}_{0} \times {\bf m}({\bf x}_{0}) + \e [{\bf x}_{0}\times {\bf m}({\bf x}_{0}))\times {\bf m}({\bf x}_{0})] = {\bf 0} $

\noindent The relation $ (a) $ can be written in the form:

\noindent $(b)~~~~~ {\bf x}_{0} \times {\bf m}({\bf x}_{0}) - \e [{\bf m}({\bf x}_{0})\times ({\bf x}_{0})\times {\bf m}({\bf x}_{0}))] = {\bf
0} $

\noindent Multiplying the relation $ (b) $ with the vector $ {\bf x}_{0}, $ we obtain:

\noindent $(c)~~~~~{\bf x}_{0}\cdot ({\bf x}_{0} \times {\bf m}({\bf x}_{0})) - \e {\bf x}_{0}\cdot [{\bf m}({\bf x}_{0})\times ({\bf
x}_{0})\times {\bf m}({\bf x}_{0}))] = {\bf 0}. $

\noindent Using the equality $(21)$, the relation $ (c) $ is equivalent with:

\noindent $(d)~~~~~ - \e ( {\bf x}_{0}\times {\bf m}({\bf x}_{0}))^{2} = {\bf 0}.$

\noindent From $ (d) $ (if $ \e \neq 0 $), follows $ {\bf x}_{0}\times {\bf m}({\bf x}_{0}) = {\bf 0}, $ that is $ {\bf x}_{0} $ is an
equilibrium point for $ (20).$\hfill\b
\bigskip

The equilibrium points of the dynamics $(12)$ are well-known (see M. Puta and D. Com\u anescu, [12]) and these are presented in the following
proposition.
\bigskip

{\bf Proposition 4.1.} $( [12] )$ {\it The equilibrium points of the Hamilton - Poisson system $ (12) $ are the following:\\
$(i)~~~~~~~e_{1} = (0,0,0);$\\[0.1cm]
$(ii)~~~~~~e_{2} = (\displaystyle\frac{a}{\lambda - a_{1}}, \displaystyle\frac{b}{\lambda - a_{2}}, \displaystyle\frac{c}{\lambda - a_{3}})~~$ for $~~\lambda\in {\bf R}\setminus \{a_{1}, a_{2}, a_{3}\}$;\\[0.2cm]
$(iii)~~~~~e_{3} = (\alpha, - \displaystyle\frac{b}{a_{2} - a_{1}}, \displaystyle\frac{c}{a_{1} - a_{3}})~~$ for $~~ \alpha \in {\bf R},~~$ if $~~ a = 0;$\\[0.2cm]
$(iv)~~~~~~e_{4} = (\displaystyle\frac{a}{a_{2} - a_{1}},\alpha ,- \displaystyle\frac{c}{a_{3} - a_{2}})~~$ for $~~ \alpha \in {\bf R},~~$ if $~~ b = 0;$\\[0.2cm]
$(v)~~~~~~~e_{5} = (- \displaystyle\frac{a}{a_{1} - a_{2}}, \displaystyle\frac{b}{a_{3} - a_{2}}, \alpha)~~$ for $~~ \alpha \in {\bf R},~~$ if $~~ c= 0. $}\\[0.2cm]

By Theorem 4.1, the equilibrium points of the $ \e$- revised system $(20)$ are $~ e_{1}, ... , e_{5}~ $ indicated in the Proposition 4.1.








It is well-known that the dynamics $(12)$ have the first integrals $ H $ and $ C $ given by $(14)$ and $(16)$. These first integrals may be
written thus:
\begin{equation}
H(x^{1}, x^{2}, x^{3}) = \displaystyle\frac{1}{2} {\bf x}\cdot {\bf I}^{-1}{\bf x} + {\bf a}\cdot {\bf x}~~~\hbox{and}~~~
C(x^{1}, x^{2}, x^{3}) = \displaystyle\frac{1}{2} {\bf x}^{2} \label{27}
\end{equation}
where ${\bf I} $ is inertia tensor and ${\bf I}^{-1} $ is its inverse. We have:
\begin{equation}
\displaystyle\frac{dH}{dt}({\bf x}) = {\bf m}({\bf x})\cdot \dot{{\bf x}}~~~\hbox{and}~~~
\displaystyle\frac{dC}{dt}({\bf x}) = {\bf x} \cdot \dot{{\bf x}}.\label{28}
\end{equation}
Indeed, $~~\displaystyle\frac{d H}{dt} = ( a_{1}x^{1} + a )\dot{x}^{1} + ( a_{2}x^{2} + b )\dot{x}^{2} + ( a_{3}x^{3} + c )\dot{x}^{3} = {\bf m}({\bf x})\cdot \dot{{\bf x}}~$ and\\[0.2cm]
$\displaystyle\frac{d C}{dt} =  x^{1}\dot{x}^{1} + x^{2}\dot{x}^{2} + x^{3} \dot{x}^{3} = {\bf x} \cdot \dot{{\bf x}}.$\\
\bigskip

{\bf Theorem 4.2.} $(i)~$ {\it For each $ \e \in {\bf R}, $ the function $ H $ given by $(14)$ is a first integral for the $ \e$- revised system $(20).$

\noindent $(ii)~~$ If $~ {\bf x} :{\bf R} \to {\bf R}^{3} $ is a solution of the $\e$- revised system, then:
\begin{equation}
\displaystyle\frac{d}{dt}(\displaystyle\frac{1}{2}{\bf x}^{2}) = - \e ( {\bf x} \times {\bf m}({\bf x}))^{2}.\label{29}
\end{equation}

\noindent $(iii)~$ For $ \e \in {\bf R}^{*}, $ the function $ C $ is not a first integral for the $ \e$- revised system.}
\bigskip

{\bf Proof.} $(i)~$ Multiplying the relation $(24)$ with the vector $ {\bf m}({\bf x}), $ we have:\\[0.1cm]
$ {\bf m}({\bf x})\cdot \dot{{\bf x}} = {\bf m}({\bf x})\cdot ({\bf x} \times {\bf m} ({\bf x})) +\e {\bf m}({\bf x})\cdot [({\bf x}\times {\bf m}({\bf x}))\times {\bf m}({\bf x})] = {\bf 0}.~$
Applying now $(28),$ we obtain
$~\displaystyle\frac{dH}{dt} = {\bf m}({\bf x})\cdot \dot{{\bf x}}= {\bf 0}. $ Hence $ H $ is a first integral for the system $(20)$.

$(ii)~$ Multiplying the relation $(24)$ with the vector $ {\bf x}, $ we have\\
$ {\bf x}\cdot \dot{{\bf x}} = {\bf x}\cdot ({\bf x} \times {\bf m} ({\bf x})) +\e {\bf x}\cdot [({\bf x}\times {\bf m}({\bf x}))\times {\bf m}({\bf x})] = -
\e {\bf x}\cdot [{\bf m}({\bf x})\times ({\bf x}\times {\bf m}({\bf x}))].$

Using now the equality $(21)$, we obtain
$~ {\bf x}\cdot \dot{{\bf x}} = - \e ({\bf x}\times {\bf m}({\bf x}))^{2}.~$
Then, we have $~ \displaystyle\frac{d}{dt}( \displaystyle\frac{1}{2}{\bf x}^{2}) =
 {\bf x}\cdot \dot{{\bf x}} = - ({\bf x}\times {\bf m}({\bf x}))^{2}.$

$(iii)~$ This assertion follows from the second relation of $(28)$ and $(ii)$.\hfill\b
\bigskip

{\bf Remark 4.1.} The function $ H $ given by $(14)$ can be put in the equivalent form:
\begin{equation}
H(x^{1}, x^{2}, x^{3}) = \displaystyle\frac{1}{2} [ a_{1}(x^{1} + \displaystyle\frac{a}{a_{1}})^{2} +  a_{2}(x^{2} + \displaystyle\frac{b}{a_{2}})^{2} + a_{3}(x^{3} + \displaystyle\frac{c}{a_{3}})^{2}] -
\displaystyle\frac{1}{2}( \displaystyle\frac{a^{2}}{a_{1}} + \displaystyle\frac{b^{2}}{a_{2}} +\displaystyle\frac{c^{2}}{a_{3}}). \label{30}
\end{equation}
For a given constant $ k\in {\bf R}, $ the geometrical image of the surface:\\[-0.2cm]
$$ H(x^{1}, x^{2}, x^{3}) = k $$\\[-0.5cm]
is an ellipsoid, since $ a_{1} > 0, a_{2} > 0, a_{3} > 0. $\hfill\b
\bigskip

{\bf Proposition 4.2.} {\it The set of equilibrium points which belong to the ellipsoid $~ H(x^{1}, x^{2}, x^{3}) = k,~ $ is finite.}
\bigskip

{\bf Proof.} Following the description of the equilibrium points given in Proposition 4.1, we remark that:

$(i)~$ the equilibrium points of the form $ ~e_{3}~ $ ( similarly, for $~ e_{4}~ $ and $~ e_{5}~ $ ) make a straight line; the intersection
between a straight line and an ellipsoid have at most two points; we deduce that  on the chosen ellipsoid there exist at most two points of the
form $~ e_{3}.$

$(ii)~$ the equilibrium points of the form $~e_{2}~$ can be obtained by solving with respect $ \lambda $ the following equation:
$$\displaystyle\frac{1}{2}[a_{1}(\displaystyle\frac{a}{\lambda -a_{1}})^{2}+ a_{2}(\displaystyle\frac{b}{\lambda
-a_{2}})^{2}+a_{3}(\displaystyle\frac{c}{\lambda -a_{3}})^{2}] + \displaystyle\frac{a^{2}}{\lambda -a_{1}} + \displaystyle\frac{b^{2}}{\lambda
-a_{2}} + \displaystyle\frac{c^{2}}{\lambda -a_{3}} = k.$$

The above equation is equivalent with the determination of roots of a polynomial of degree at most $ 6; $ therefore on the chosen ellipsoid
there exist at most $ 6 $ equilibrium points of the form $~ e_{2}.$\hfill\b

\section {The behaviour of the solutions of the $ \e $ - revised system }
\textbf{Theorem 5.1} (i) The solutions of the $\varepsilon$-revised system are bounded.

\noindent (ii) The maximal solutions of the $\varepsilon$-revised system are globally solutions (i.e. these are defined on $\mathbf{R}$).
\bigskip

\textbf{Proof} (i) Given a solution of (20), there exists a constant $k$ such that its trajectory lie on the ellipsoid $H(x^{1},x^{2},x^{3})=k$.
From this we deduce that all solutions are bounded.

(ii) Let $\textbf{x}:(m,M)\subset \mathbf{R}\rightarrow \mathbf{R}^{3}$ be a maximal solution. We assume that $\textbf{x}$ is not globally. It
follows $m>-\infty$ or $M<\infty$. In these situations, we known that there exists $k\in \mathbf{R}$ such that $H(x^{1},x^{2},x^{3})=k$ for all
$t\in \mathbf{R}$ and the graph of the solution is contained in a compact domain. According with [6] (theorem 3.2.5, p.141) we obtain a
contradiction with the fact that $\textbf{x}$ admit a prolongation on the right or the left (also, can be applied the theorem of Chilingworth
(1976), see theorem 1.0.3, p.7 in [3]). \hfill\b
\bigskip

In the sequel we study the asymptotic behaviour of the globally solutions of the $\varepsilon$-revised system.

Denote by $\mathbf{E}$ the set of equilibrium points of the $\varepsilon$-revised system (20) and by $\Gamma$ the trajectory of a solution
$\textbf{x}:\mathbf{R}\rightarrow\mathbf{R}^{3}$ of (20). By theory of differential equations (see [10] p. 174-176), the $\omega$-limit set and
$\alpha$-limit set of $\Gamma$ are:
$$\omega (\Gamma)=\{\textbf{y}\in \mathbf{R}^{3}\,\,/\,\,\exists t_{n}\rightarrow\infty\,\,\texttt{such
that}\,\,\textbf{x}(t_{n})\rightarrow\textbf{y}\},$$
$$\alpha (\Gamma)=\{\textbf{z}\in \mathbf{R}^{3}\,\,/\,\,\exists t_{n}\rightarrow -\infty\,\,\texttt{such
that}\,\,\textbf{x}(t_{n})\rightarrow\textbf{z}\}.$$
\smallskip

\textbf{Theorem 5.2} Let $\textbf{x}:\mathbf{R}\rightarrow\mathbf{R}^{3}$ be a solution of the $\varepsilon$-revised system with
$\varepsilon\neq 0$. There exist the equilibrium points $\textbf{x}_{m},\textbf{x}_{M}\in \mathbf{R}^{3}$ of the system (20) such that
$\lim_{t\rightarrow -\infty}\textbf{x}(t)=\textbf{x}_{M}$ and $\lim_{t\rightarrow \infty}\textbf{x}(t)=\textbf{x}_{m}$.
\bigskip

\textbf{Proof} The theorem is proved in the following steps:

\noindent (i) $\alpha (\Gamma)\neq \emptyset$ and $\omega (\Gamma)\neq \emptyset$.

\noindent (ii) $\alpha (\Gamma)\bigcap \omega (\Gamma)\subset \mathbf{E}$.

\noindent (iii) The sets $\alpha (\Gamma)$ and $\omega (\Gamma)$ contains exactly one element.

Taking account into that each solution is bounded (hence it is contained in a compact domain) and applying theorem 1, p. 175 in [10], we obtain
immediately the assertions (i).

(ii) For demonstration consider the case when $\varepsilon >0$. Using the relation (29), we deduce that the function $t\rightarrow
\textbf{x}^{2}(t)$ is a strictly decreasing function. Being bounded it follows that there exists $\lim_{t\rightarrow\infty}\textbf{x}^{2}(t)=L$
and $L$ is finite.

For each $\textbf{y}\in \omega (\Gamma)$ there exists the sequence $t_{n}\rightarrow\infty$ such that $\textbf{x}(t_{n})\rightarrow\textbf{y}$.
Then $\textbf{x}^{2}(t_{n})\rightarrow\textbf{y}^{2}$ and hence $\textbf{y}^{2}=L$.

By theorem 2, p.176 in [10], we have that the trajectory $\Gamma_{\textbf{y}}$ of the solution $\textbf{x}_{\textbf{y}}$ which verifies the
initial condition $\textbf{x}_{\textbf{y}}(0)=\textbf{y}$, satisfies the relation $\Gamma_{\textbf{y}}\subset \omega (\Gamma)$.

If we assume that $\textbf{y}$ is not an equilibrium point, then we deduce (using the relation (29)) that for $t>0$ we have
$\textbf{x}_{\textbf{y}}^{2}(t)<L$ and this is in contradiction with the above result. Therefore, we have $\omega (\Gamma)\subset \mathbf{E}$.

Similarly, we prove that $\alpha (\Gamma)\subset \mathbf{E}$. Hence the assertion (ii) holds.

The case $\varepsilon<0$ is similar.

(iii) There exists a constant $k$ such that the sets $\alpha (\Gamma)$ and $\omega (\Gamma)$ are included in the ellipsoid
$H(x^{1},x^{2},x^{3})=k$. By (ii), we deduce that $\alpha (\Gamma)$ and $\omega (\Gamma)$ are included in the set of equilibrium points which
lies of the above ellipsoid. On the other hand, applying Proposition 4.2 and using the fact that $\alpha (\Gamma)$ and $\omega (\Gamma)$ are
connected (see theorem 1, p.175 in [10]), we obtain that $\alpha (\Gamma)$ and $\omega (\Gamma)$ are formed by only one element.\hfill\b
\bigskip

\textbf{Remark 5.1} Using the relation (29) it is easy to observe that the following assertions hold:

(i) if $\varepsilon >0$ then $\textbf{x}^{2}_{M}>\textbf{x}^{2}_{m}$;

(ii) if $\varepsilon <0$ then $\textbf{x}^{2}_{M}<\textbf{x}^{2}_{m}$.\hfill\b
\bigskip

As an immediate consequence we obtain the following theorem.
\bigskip

\textbf{Theorem 5.3} If $\varepsilon\neq 0$, then for each solution $\textbf{x}:\mathbf{R}\rightarrow\mathbf{R}^{3}$ of the
$\varepsilon$-revised system we have:
\begin{equation}\label{*}
    \left\{%
\begin{array}{ll}
    \texttt{if}\,\, t\rightarrow \infty \,\,\Rightarrow\,\, d(\textbf{x}(t),\mathbf{E})\rightarrow 0 \\
    \texttt{if}\,\, t\rightarrow -\infty \,\,\Rightarrow\,\, d(\textbf{x}(t),\mathbf{E})\rightarrow 0. \\
\end{array}%
\right.
\end{equation}\hfill\b
\smallskip

\textbf{Remark 5.2} From Theorem 5.3 follows that the set $\mathbf{E}$ of the equilibrium points is an attracting set (see definition 2, p.178
in [10]) and also is a reppeling set (see [3], p.34). Thus, the space $\mathbf{R}^{3}$ is simultaneously a domain of attraction and a domain of
repulsion of $\mathbf{E}$. \hfill\b

\section {The Lyapunov stability of equilibrium points of the $ \e $ - revised system in the case $\e > 0 $}

{\bf The stability of the point $ e_{1} = ( 0, 0, 0 ) $}. We have the following results.
\bigskip

\textbf{Theorem 6.1} The equilibrium point $ e_{1}$ is Lyapunov stable.
\bigskip

\textbf{Proof} Let $\gamma >0,\, t_{0}\in \mathbf{R}$ and $\textbf{x}_{0}\in \mathbf{R}^{3}$ such that $|\textbf{x}_{0}|<\gamma$, where $|\cdot
|$ denotes the euclidian norm in $\mathbf{R}^{3}$. Denote by $t\rightarrow \textbf{x}(t,t_{0},\textbf{x}_{0})$ the solution of the
$\varepsilon$-revised system which verifies the initial condition $\textbf{x}(0,t_{0},\textbf{x}_{0})=\textbf{x}_{0}$.

Using the relation $|\textbf{x}(t,t_{0},\textbf{x}_{0})|=\sqrt{\textbf{x}^{2}(t,t_{0},\textbf{x}_{0})}$ and according with the relation (29), we
observe that the function $t\rightarrow \textbf{x}(t,t_{0},\textbf{x}_{0})$ is a decreasing function and hence we have:
$$|\textbf{x}(t,t_{0},\textbf{x}_{0})|\leq |\textbf{x}_{0}|<\gamma\,\, \texttt{for}\,\, t>t_{0}.$$
Then (see [4], p.22) we have that $e_{1}$ is a Lyapunov stable equilibrium point.\hfill\b
\bigskip

\textbf{Remark 6.1} The equilibrium point $e_{1}$ is not asymptotical stable.

Indeed, if $a=b=c=0$ then the coordinates axis are formed from equilibrium points. If al least of one of the numbers $a,b,c$ is non null, then:
$$\texttt{if}\,|\lambda|\rightarrow \infty\,\Rightarrow\,(\frac{a}{\lambda -a_{1}}, \frac{b}{\lambda -a_{2}},\frac{c}{\lambda
-a_{3}})\rightarrow (0,0,0).$$ Hence, in all neighbourhood of $e_{1}$ there exist an infinity of equilibrium points.\hfill\b
\bigskip

{\bf The stability of the point $\overline{\textbf{x}_{0}}= ( - \displaystyle\frac{a}{a_{1}}, - \displaystyle\frac{b}{a_{2}}, -
\displaystyle\frac{c}{a_{3}} ) $}. The equilibrium point $\overline{\textbf{x}_{0}}$ is an equilibrium point of the form $e_{2}$ and it is
obtained for $\lambda =0$.
\bigskip

\textbf{Theorem 6.2} The equilibrium point $\overline{\textbf{x}_{0}}$ is Lyapunov stable.
\bigskip

\textbf{Proof} Using the relation (30) and the inequality $0<a_{1}<a_{2}<a_{3}$, we deduce:
$$\frac{a_{1}}{2}|\textbf{x}-\overline{\textbf{x}_{0}}|\leq H(\textbf{x})-H(\overline{\textbf{x}_{0}}) \leq \frac{a_{3}}{2}|\textbf{x}-\overline{\textbf{x}_{0}}|$$
For $t_{0}\in \mathbf{R}$ and $\textbf{x}_{0}\in \mathbf{R}^{3}$ denote with $\textbf{x}(t,t_{0},\textbf{x}_{0})$ a solution of
$\varepsilon$-revised system which verifies the initial condition $\textbf{x}(0,t_{0},\textbf{x}_{0})=\textbf{x}_{0}$.

Let $\gamma >0$ and $\delta (\gamma)=\frac{2\gamma}{a_{1}}$. Let $\textbf{x}_{0}\in \mathbf{R}^{3}$ such that:
$$ H(\textbf{x}_{0})-H(\overline{\textbf{x}_{0}})\leq \delta (\gamma).$$
From the fact that $H$ is a first integral we deduce that:
$$H(\textbf{x}(t,t_{0},\textbf{x}_{0})-H(\overline{\textbf{x}_{0}})= H(\textbf{x}_{0})-H(\overline{\textbf{x}_{0}}).$$
Hence for all $t\in \mathbf{R}$ the following inequality holds:
$$\frac{a_{1}}{2}|\textbf{x}-\overline{\textbf{x}_{0}}|\leq \delta (\gamma)$$
and we obtain that $\overline{\textbf{x}_{0}}$ is Lyapunov stable. \hfill\b
\bigskip

\textbf{Remark 6.2} The stable equilibrium point $\overline{\textbf{x}_{0}}$ realizes the absolute minimum of the function $H$. \hfill\b
\bigskip

{\bf The unstability of equilibrium points of the form $e_{2}$  with $\lambda \in (0, a_{1})$}. For the demonstration of this results we use the
Theorem 6.3 and Lemma 6.1.
\bigskip

\textbf{Theorem 6.3} If $\textbf{x}_{0}\in \mathbf{E}$ such that there exists $\textbf{y}\in \mathbf{E}$ with the properties:
 \begin{center}
 (i) $H(\textbf{y})=H(\textbf{x}_{0})$ and (ii) $|\textbf{y}|<|\textbf{x}_{0}|$\end{center}
then $\textbf{x}_{0}$ is an unstable equilibrium point.
\bigskip

\textbf{Proof} For $k\in \mathbf{R}$ denote by $\mathbf{E}_{k}=\{\textbf{x}\in \mathbf{E}\,/\,H(x)=k\}$. The set
$\mathbf{E}_{H(\textbf{x}_{0})}$ is finite (by Proposition 4.2). We denote:
$$\gamma_{0} =\min\{|\textbf{x}-\textbf{x}_{0}| / \textbf{x}\in E_{H(\textbf{x}_{0})}-\{\textbf{x}_{0}\}\}$$
Let $\textbf{z}\in \mathbf{R}^{3}$ such that $H(\textbf{z})=H(\textbf{x}_{0})$ and $|\textbf{z}|<|\textbf{x}_{0}|$. Then:
$$\lim_{t\rightarrow \infty}\textbf{x}(t,0,\textbf{z})\in \mathbf{E}$$
and
$$\texttt{if} \,\,t>0 \,\Rightarrow\, |\textbf{x}(t,0,\textbf{z})|<|\textbf{z}|$$
and we deduce that there exists $t_{z}>0$ such that:
$$|\textbf{x}(t,0,\textbf{z})-\textbf{x}_{0}|>\frac{\gamma_{0}}{2}\,\, \texttt{if}\,\,t>t_{z}$$
It follows that $\textbf{x}_{0}$ is unstable. \hfill\b
\bigskip

We assume that $(a,b,c)\neq (0,0,0)$ and we introduce the notation:
$$e_{2\lambda}=(\frac{a}{\lambda -a_{1}},\frac{b}{\lambda -a_{2}},\frac{c}{\lambda -a_{3}})\,\,\texttt{for all}\,\,\lambda\in \mathbf{R}-\{a_{1},a_{2},a_{3}\}$$
\bigskip

\textbf{Lemma 6.1} (i) If $\sigma<\mu<a_{1}$, then $|e_{2\sigma}|<|e_{2\mu}|$.

\noindent (ii) If $\sigma,\mu>0,\,\,(\frac{\mu}{\sigma})^{2}>\frac{a_{3}}{a_{1}}$ and $H(e_{2\sigma})=H(e_{2\mu})$, then
$|e_{2\sigma}|>|e_{2\mu}|$.

\noindent (iii) If $0<\sigma<a_{1}<a_{3}<\mu$ and $H(e_{2\sigma})=H(e_{2\mu})$, then $|e_{2\sigma}|>|e_{2\mu}|$.
\bigskip

\textbf{Proof} (i) Consider the function $g:(-\infty, a_{1})\rightarrow \mathbf{R}$ given by:
$$g(\lambda)=(\frac{a}{\lambda -a_{1}})^{2}+(\frac{b}{\lambda -a_{2}})^{2}+(\frac{c}{\lambda -a_{3}})^{2}.$$
The derivative of the function $g$ is:
$$g'(\lambda)=-\frac{2a^{2}}{(\lambda -a_{1})^{3}}-\frac{2b^{2}}{(\lambda -a_{2})^{3}}-\frac{2c^{2}}{(\lambda -a_{3})^{3}}$$
We observe that $g'(\lambda)>0$ and we obtain that $g$ is a strictly increasing function. We have:
$$g(\sigma)=|e_{2\sigma}|^{2},\,\,g(\mu)=|e_{2\mu}|^{2}$$
and we obtain the desired result.

(ii) From hypothesis $H(e_{2\sigma})=H(e_{2\mu})$ follows that there exists a constant $q>0$ with the following properties:
$$\frac{1}{a_{1}}\frac{a^{2}}{(\sigma -a_{1})^{2}}+\frac{1}{a_{2}}\frac{b^{2}}{(\sigma -a_{2})^{2}}+\frac{1}{a_{3}}\frac{c^{2}}{(\sigma
-a_{3})^{2}}=\frac{q}{\sigma ^{2}}$$
$$\frac{1}{a_{1}}\frac{a^{2}}{(\mu -a_{1})^{2}}+\frac{1}{a_{2}}\frac{b^{2}}{(\mu -a_{2})^{2}}+\frac{1}{a_{3}}\frac{c^{2}}{(\mu
-a_{3})^{2}}=\frac{q}{\mu ^{2}}$$ Using $a_{1}<a_{2}<a_{3}$, we obtain the inequalities:
$$|e_{2\sigma}|^{2}>\frac{a_{1}q}{\sigma},\,\,|e_{2\mu}|^{2}<\frac{a_{3}q}{\mu}$$
and we observe that the assertion (ii) holds.

(iii) This assertion follows immediately from (ii).\hfill\b
\bigskip

\textbf{Theorem 6.4} The equilibrium point $e_{2\lambda}$ with $0<\lambda <a_{1}$ is unstable.
\bigskip

\textbf{Proof} Consider the function $h:(-\infty,a_{1})\bigcup(a_{3},\infty)\rightarrow \mathbf{R}$ given by:
$$h(\sigma)=H(e_{2\sigma}).$$
Using the relation (30) for $H$, we find:
$$h(\sigma)=\frac{\sigma^{2}}{2}[\frac{a^{2}}{a_{1}(\sigma-a_{1})^{2}}+\frac{b^{2}}{a_{2}(\sigma-a_{2})^{2}}+\frac{c^{2}}{a_{3}(\sigma-a_{3})^{2}}]-\frac{1}{2}(\frac{a^{2}}{a_{1}}+\frac{b^{2}}{a_{2}}+\frac{c^{2}}{a_{3}})$$
The function $h$ have the following properties:
\begin{itemize}
    \item 0 is an absolute minimum point.
    \item $\lim_{\sigma\rightarrow -\infty}h(\sigma)=\lim_{\sigma\rightarrow \infty}h(\sigma)=0$.
    \item $\lim_{\sigma\rightarrow a_{1}}h(\sigma)=\lim_{\sigma\rightarrow a_{3}}h(\sigma)=\infty$.
    \item $h$ is strictly decreasing on $(-\infty, 0)$, strictly increasing on $(0,a_{1})$ and strictly decreasing on $(a_{3},\infty)$.
\end{itemize}

The demonstrations divided on three cases.

(I) Assume that $h(\lambda)<0$. In this situation there exists $\sigma<0<\lambda<a_{1}$ such that $h(\lambda)=h(\sigma)$ and imply
$H(e_{2\lambda})=H(e_{2\sigma})$. Hence the equilibrium points $e_{2\lambda}$ and $e_{2\sigma}$ belong to same ellipsoid.

On the other hand, by Lemma 6.1 (ii), follows $|e_{2\sigma}|<|e_{2\lambda}|$. Applying now Theorem 6.3, deduce that $e_{2\lambda}$ is an
unstable equilibrium point.

(II) Assume that $h(\lambda)=0$ we have $H(e_{2\lambda})=H(0,0,0)$ and it is clearly that $|(0,0,0)|<|e_{2\lambda}|$. By Theorem 6.3 we find the
desired result.

(III) Assume that $h(\lambda)>0$. Then there exists $\sigma>a_{3}$ such that $h(\lambda)=h(\sigma)$ and hence $H(e_{2\lambda})=H(e_{2\sigma})$.
Applying Lemma 6.1 (iii) follows $|e_{2\lambda}|>|e_{2\sigma}|$ and by Theorem 6.3 we deduce that $e_{2\lambda}$ is unstable.\hfill\b
\bigskip

{\bf The stability of equilibrium points of the form $ e_{2} $  with $ \lambda < 0 $. }
\bigskip

\textbf{Theorem 6.5} The equilibrium points of the form $e_{2}$ with $\lambda<0$ are Lyapunov stables.
\bigskip

\textbf{Proof} Let $\lambda<0$ and the equilibrium point
$\textbf{x}_{0}=(\frac{a}{\lambda-a_{1}},\frac{b}{\lambda-a_{2}},\frac{c}{\lambda-a_{3}})$ of the form $e_{2}$. It is well-known that the study
of stability of $\textbf{x}_{0}$ in the Lyapunov sense is equivalent with the study of stability of the null solution $(0,0,0)$ for the
differential system obtained from the $\varepsilon$-revised system by transformation of variables:
$$\textbf{z}=\textbf{x}-\textbf{x}_{0}$$
The system obtained in this manner is called the \textbf{perturbed $\varepsilon$-revised system}.

Consider the function $K:\mathbf{R}^{3}\rightarrow\mathbf{R}$ given by
$$K(\textbf{z})=\frac{1}{2}\textbf{z}\cdot\textbf{I}^{-1}\textbf{z}-\frac{\lambda}{2}\textbf{z}^{2}$$
Since the tensor $\textbf{I}^{-1}$ is strictly positive definite and $\lambda<0$ we obtain that $K$ is a quadratic form strictly positive
definite.

Next we prove that if $\textbf{z}:\mathbf{R}\rightarrow\mathbf{R}^{3}$ is a solution for the perturbed $\varepsilon$-revised system, then:
$$\frac{d}{dt}K(\textbf{z}(t))<0$$
By a direct computation and taking account into the relations (27) we have:
$$K(\textbf{z})=H(\textbf{x})-\lambda
C(\textbf{x})-\frac{1}{2}\textbf{x}_{0}\cdot\textbf{I}^{-1}\textbf{x}_{0}-\textbf{a}\cdot\textbf{x}_{0}+\frac{\lambda}{2}\textbf{x}_{0}^{2}$$
Applying now Theorem 4.2, we obtain:
$$\frac{d}{dt}K(\textbf{z}(t))=\varepsilon\lambda(\textbf{x}\times\textbf{m(x)})^{2}$$
and follows that $\frac{d}{dt}K(\textbf{z}(t))<0$, since $\varepsilon>0,\,\,\lambda<0$.

It is easy to see that $K^{*}(t)=K(\textbf{z}(t))$ is a strictly decreasing function. By theorem 1.1, p.21 in the paper [4] we deduce that
$\textbf{x}_{0}$ is Lyapunov stable.\hfill\b
\bigskip

{\bf Conclusion - the stability problem for the Hamilton Poisson system $(12)$ versus the $ \e $- revised system $(20) $ with $\e > 0 $}

Concerning to the equilibrium points of the system $(12)$ are established the following results (see, theorem 1.1, [12]):\\
$(1)~~~~~e_{1}~$ is Lyapunov stable;\\
$(2)~~~~~e_{2}~$ are Lyapunov stables for $ \lambda \in (- \infty, a_{1} ) \cup ( a_{3}, \infty )$;\\
$(3)~~~~~e_{3}~$ are Lyapunov stables;\\
$(4)~~~~~e_{4}~$ are unstables;\\
$(5)~~~~~e_{5}~$ are Lyapunov stables.

By Remark 4.1 (see, [12]), there exist cases for which the problem to decide the nonlinear stability or unstability are not discussed.

For the stability of equilibrium points of the $ \e$- revised system $(20)$ with $ \e > 0 $ have proved the following assertions:\\
$(1)~~~~~e_{1}~$ is Lyapunov stable;\\
$(2)~~~~$ the equilibrium points of the form $~e_{2}~$ with $ \lambda \leq 0 $ are Lyapunov stables;\\
$(3)~~~~$ the equilibrium points of the form $~e_{2}~$ with $ 0 < \lambda < a_{1}~ $ are unstables.
\bigskip

\begin{center}
{\bf References}\\
\end{center}
\noindent
\hspace*{0.7cm} {\bf [1]. R. Abraham, J.E. Marsden}, {\it Foundations of Mechanics.Second Edition}. Addison-Wesley, 1978.\\
\hspace*{0.7cm} {\bf [2]. D. Fish},{\it Dissipative perturbation of 3D Hamiltonian systems.Metriplectic systems}. Preprint, arXiv:math-ph/0506047, v1,2005.\\
\hspace*{0.7cm} {\bf [3]. J. Guckenheimer, P. Holmes}, {\it Nonlinear oscilations, dynamical systems and bifurcations of vector fields}. Springer- Verlag, New York, 1990.\\
\hspace*{0.7cm} {\bf [4]. A. Halanay}, {\it Teoria calitativ\u a a ecua\c tiilor diferen\c tiale }. Ed. Academiei, Bucure\c sti, 1963.\\
\hspace*{0.7cm} {\bf [5] Gh. Ivan, D. Opri\c s}, {\it Dynamical systems on Leibniz algebroids}. Differential Geometry - Dynamical Systems}, {\bf 8}( 2006 ), 127 - 137.\\
\hspace*{0.7cm} {\bf [6]. St. Miric\u a}, {\it Ecua\c tii diferen\c tiale \c si cu derivate par\c tiale I}. Litografia Univ. Bucure\c sti, 1989.\\
\hspace*{0.7cm} {\bf [7] J.E. Marsden}, {\it Lectures on Mechanics}. London Mathematical Society, Lectures Note Series, vol.{\bf 174}, 2 nd edition, Cambridge University Press, 1992.\\
\hspace*{0.7cm} {\bf [8] P.J. Morrison}, {\it A paradigm for joined Hamiltonian and dissipative systems}. Physica,{\bf 18}D ( 1986), 410 - 419.\\
\hspace*{0.7cm} {\bf [9]. J.- P. Ortega, V. Planas -Bielsa}, {\it Dynamics on Leibniz manifolds}. Preprint, arXiv:math. DS/0309263,2003.\\
\hspace*{0.7cm} {\bf [10]. L. Perko}, {\it Differential equations and dynamical systems}. Springer- Verlag, New York, 1991.\\
\hspace*{0.7cm} {\bf [11]. M. Puta}, {\it Hamiltonian mechanics and geometric quantization}. Mathematics and its Applications, vol. {\bf 260}, Kluwer, 1993.\\
\hspace*{0.7cm} {\bf [12]. M. Puta, D. Com\u anescu}, {\it On the rigid body with three linear controls}. Analele Univ. din Timi\c soara, Seria Matematic\u a - Informatic\u a, vol. {\bf 35} (1), 1997, p.-.\\
\hspace*{0.7cm} {\bf [13]. M. Puta, D. Com\u anescu, S. Chirici}, {\it Elemente de mecanic\u a hamiltonian\u a}. Ed. Mirton, Timi\c soara, 2004.\\[0.2cm]

\hspace*{0.7cm} Seminarul de Geometrie \c si Topologie\\
\hspace*{0.7cm} West University of Timi\c soara\\
\hspace*{0.7cm} Bd-ul V. P{\^a}rvan no.4, 300223, Timi\c soara\\
\hspace*{0.7cm} Romania\\
\hspace*{0.7cm} E-mail: comanescu@math.uvt.ro,  ivan@math.uvt.ro  and  mihai31ro@yahoo.com\\[0.2cm]

\end{document}